\documentclass{amsart}
\usepackage{array}

\def\C{\mathbb C}
\def\P{\mathbb P}

\def\Z{\mathbb Z}
\def\Q{\mathbb Q}

\def\F{\mathbb F}

\def\bmu{\mathbf{\mu}}

\newtheorem{Theorem}{Theorem}
\newtheorem{Lemma}{Lemma}
\newtheorem{Corollary}{Corollary}
\newtheorem{Remark}{Remark}
\newtheorem{Example}{Example}
\newtheorem{Definition}{Definition}

\def\qed{\hfill \ensuremath{\Box}}
\def\Gal{\mathrm{Gal\,}}
\def\Qbar{\overline{\Q}}

\begin{document}
\title{The modularity of K3 surfaces with non-symplectic group actions}
%[The modularity of K3 surfaces] 
\author{Ron Livn\'e}
\address{Institute of Mathematics \\
  Hebrew University of Jerusalem Givat-Ram\\
   Jerusalem 91904\\
  Israel}
\email{rlivne@math.huji.ac.il}

\author{Matthias Sch\"utt}
\address{Institute for Algebraic Geometry\\
 Leibniz University Hannover\\
 Welfengarten 1\\
 30167 Hannover\\
 Germany}
\email{schuett@math.uni-hannover.de}
%\urladdr{http://www.iag.uni-hannover.de/\~{}schuett/}

\author{Noriko Yui}
\address{Department of Mathematics and Statistics\\
  Queen's University\\
  Kingston, Ontario K3L 3N6\\
  Canada}
\email{yui@mast.queensu.ca}
\date{\today}
\thanks{R.~Livn\'e was partially supported by an ISF grant.
M.~Sch\"utt was supported by Deutsche Forschungsgemeinschaft (DFG) under grant Schu 2266/2-2. 
  N.~Yui was supported in part by
  Discovery Grant of the Natural Sciences and Engineering Research 
  Council (NSERC) of Canada}

\subjclass[2000]{Primary 14J32, 14J28, 14J27, 14J20, 14G10, 11G40,11F80} 

\keywords{K3 surface, non-symplectic group action,
 automorphism group of a K3 surface, unimodular and non-unimodular lattices,
Galois representation, modularity} 

\begin{abstract}
We consider complex K3 surfaces with a non-symplectic group acting trivially on the algebraic cycles.
Vorontsov and Kond\=o classified those K3 surfaces with 
transcendental lattice of minimal rank.
The purpose of this note is to study the Galois representations associated
to these K3 surfaces. The rank of the transcendental
lattices is even and varies from $2$ to $20$, excluding $8$ and $14$. We show that these K3 
surfaces are dominated by Fermat surfaces, and hence they are all of CM type. We will
establish the modularity of the Galois representations associated to them.
Also we discuss mirror symmetry for these K3 surfaces in the sense of Dolgachev, and
show that a mirror K3 surface exists with one exception.
\end{abstract}
\maketitle

\begin{section}{Introduction} 
Let $X$ be an algebraic K3 surface, and let $S_X$ and $T_X$ be the lattices of algebraic and
transcendental cycles on $X$, respectively.  The automorphism group $\mbox{Aut}(X)$
acts on these lattices. Let $O(S_X)$ denote the group of isometries of $S_X$.
Then the kernel $H_X:=\mbox{Ker}(\mbox{Aut}(X)\to O(S_X))$ 
is a finite cyclic group.  
If $g$ is a generator of $H_X$, and if $\omega$ 
is a nowhere vanishing holomorphic $2$-form on $X$, then $g$ acts on $\omega$ by multiplication by a primitive
$k$-th root of unity for some $k$.  Let $\phi(k)$ denote the Euler $\phi$-function. 
Then $\phi(k)$ divides $\mbox{rank}(T_X)$. 
We are interested in the special situation where $\mbox{rank}(T_X)=\phi(k)$. There are two cases (cf.~Thm.~1, 2): 
\smallskip

(I)  $T_X$ is unimodular (i.e, $\mbox{det }(T_X)=\pm 1$). 
Then there are exactly six values for $k$, namely, $k\in\{66, 44, 42, 36, 28, 12\}$. 
Conversely, for each $k$, there exists a unique K3 surface $X$ up to isomorphism with these properties. These K3 surfaces can be defined over $\Q$ by Weierstrass equations. %over $\Q$
%The rank of $T_X$ is $20,20,12,12,12,4$, respetively.
\smallskip

(II) $T_X$ is non-unimodular. Then
there are exactly ten values for $k$, namely, $k\in\{3,9,27,5,25,7,11,13,17,19\}$.
Conversely, for each $k$, there exists a unique K3 surface $X$ up to isomorphism
with these properties. With one exception, these K3 surfaces are defined by
Weierstrass equations over $\Q$. % over $\Q$ with one exception of $k=25$. 
For $k=25$, the K3 surface
is defined as a double sextic over $\Q$, i.e.~a double cover of $\P^2$ branched along a sextic curve.

Our goal of this paper is to establish the modularity of the Galois representations 
associated to the transcendental lattices of these K3 surfaces.  We note that the 
rank of $T_X$ takes even values from $2$ to $20$, excluding the
values $8$ and $14$.  We show in Corollary \ref{Cor:CM} that the K3 surfaces are all of CM type. 
Indeed they are Delsarte surfaces. As such, they are realized as quotients of Fermat surfaces 
of appropriate degrees (Lemma~\ref{Lem:G}). We verify the modularity explicitly in 
Lemma \ref{Lem:R_t}. We study the Galois representations associated
to the transcendental lattices of our K3 surfaces over $\Q$ (Lemma 4).
 The paper concludes with a discussion of mirror symmetry for these
K3 surfaces in the sense of Dolgachev \cite{Dol}, and its arithmetic aspects.
All K3 surfaces with one exception ($k=3$) admit mirror partners (Lemma 5).
When $k\neq 3,7,11,17,19$, there are mirrors that are again Delsarte surfaces (Lemma 6).

\medskip

K3 surfaces with non-symplectic automorphisms of $2$-power order are classified in a recent article  by the second author \cite{Sch2}. The paper also considers analogous arithmetic questions.
\end{section}

\begin{section}{Classification of K3 surfaces with non-sympletic group actions} 
\label{s:results}

Let $X$ be an algebraic K3 surface over $\C$. Let $H^0(X,\Omega_X^2)$
be the space of global holomorphic $2$-forms on $X$, which is of complex
dimension $1$, i.e., $H^0(X,\Omega_X^2)\simeq\C\omega$. Let $S_X$ and $T_X$ be the lattices
of algebraic and transcendental cycles of $X$, repectively.  Let $\mbox{Aut}(X)$
be the automorphism group of $X$. Since any automorphism $g\in \mbox{Aut}(X)$ preserves
the one-dimensional vector space $H^0(X,\Omega_X^2)=\C\omega$, there is a non-zero complex number $\alpha(g)\in \C^*$
such that $g^*(\omega)=\alpha(g)\omega$. Nikulin proved that $\alpha(\mbox{Aut}(X))$ is a finite cyclic group \cite[Thm.~3.1]{Ni}. Moreover, if $k=|\alpha(\mbox{Aut}(X))|$, then $\phi(k) \,|\,\mbox{rank}(T_X)$.

\medskip

Let $O(L)$ denote the group of isometries of a lattice $L$. 
It is a consequence of the Torelli theorem that the natural representation
$$\rho\, :\, \mbox{Aut}(X)\to O(H^2(X,\Z))$$
is faithful, i.e.~the induced map $\mbox{Aut}(X)\to O(S_X)\times O(T_X)$ is injective.
We are interested in automorphisms which act trivially on $S_X$. We define the subgroup $H_X\subseteq \mbox{Aut}(X)$ by the following short exact sequence
\begin{equation*}
1\to H_X \to \mbox{Aut}(X)\to O(S_X)\to 1.
\end{equation*}
Since $\rho$ is faithful, $H_X$ can be identified with its image under $\alpha$, a finite cyclic group. 
Hence if $k=|H_X|$, then $\phi(k)$ divides the rank of $T_X$. In the following, we let always $g_k$ denote a generator for
the cyclic group $H_X$.

\medskip

There are only finitely many values for $k$ arising this way. The classification of all possible choices for such $k$ was announced by Vorontsov \cite{Vo}. It was proved later  by Vorontsov using the theory of cyclotomic fields over $\Q$, 
and by Kond\=o \cite{Ko} using the theory of elliptic surfaces
and results of Nikulin \cite{Ni} on finite automorphism groups of K3 surfaces.

\begin{Theorem} [Unimodular cases]
{\sl Let $X$ be an algebraic K3 surface over $\C$. 

Assume that $T_X$ is unimodular (i.e., $\mbox{det}(T_X)=\pm 1$). Denote $k=|H_X|$. 
Let $\Sigma=\{66,44,42,36,28,12\}$. 
Then the following assertions hold.

(a) $k$ is a divisor of an element in
$\Sigma$. In particular $k\leq 66$.

(b)  Furthermore, assume that $\phi(k)=\mbox{rank}(T_X)$.  Then $k\in\Sigma$. 

(c) Conversely, for each $k\in\Sigma$,
there exists a unique K3 surface $X$ (up to isomorphism) satisfying
the properties that $|H_X|=k$ and $\mbox{rank}(T_X)=\phi(k)$. 
For these values of $k$, $S_X$ and $T_X$ are given as follows.}
\end{Theorem}

\begin{table}[ht!]
\[
\begin{array}{c||c|c||c|c} \hline
k & S_X & \mbox{rank}(S_X) & T_X& \mbox{rank}(T_X) \\ \hline\hline
66 & U_2 & 2 & U_2^2\oplus(-E_8)^2 & 20 \\ \hline
44 & U_2 & 2 & U_2^2\oplus(-E_8)^2 & 20 \\ \hline
42 & U_2\oplus (-E_8) & 10 & U_2^2\oplus(-E_8) & 12 \\ \hline
36 & U_2\oplus(-E_8) & 10 & U_2^2\oplus(-E_8) & 12 \\ \hline
28 & U_2\oplus(-E_8) & 10 & U_2^2\oplus(-E_8) & 12 \\ \hline
12 & U_2\oplus(-E_8)^2 & 18 & U_2^2 & 4 \\
\hline
\end{array}
\]
\caption{Unimodular non-symplectic K3 surfaces}
\end{table}
\medskip

Here $U_2$ denotes the hyperbolic plane and $E_8$ is the positive-definite unimodular root lattice of rank 8. In this notation, the K3 lattice is $H^2(X,\Z)=U_2^3 \oplus (-E_8)^2$.

\medskip

The K3 surfaces in Theorem 1 can all be defined over $\Q$. Kond\=o \cite{Ko} %(cf. Machida and Oguiso \cite{MO}) 
exhibits defining equations in terms of Weierstrass models. We reproduce them in the following table up to some signs. We also list the action of a generator $g_k$ of $H_X$ on $X$. In terms of the local coordinates of the Weierstrass model, the action of $g_k$  is given by
$$(x,y,t)\mapsto (\zeta_k^{\alpha} x, \zeta_k^{\beta} y, \zeta_k^{\gamma} t)$$ 
where $\zeta_k$ is a primitive $k$-th root of unity, and
$\alpha, \beta, \gamma\in \Z/k\Z$.

\begin{table}[ht!]
\[
\begin{array}{c|c|c} \hline
k & X & g_k \\ \hline\hline
66 & y^2=x^3-t(t^{11}+1) & (x,y,t)\mapsto (\zeta_{66}^2x, \zeta_{66}^3 y, \zeta_{66}^6 t) \\ \hline
44 & y^2=x^3+x+t^{11}    & (x,y,t)\mapsto (-x, \zeta_{44}^{11}y, \zeta_{44}^2t) \\ \hline
42 & y^2=x^3-t^5(t^7+1)  & (x,y,t)\mapsto (\zeta_{42}^2x, \zeta_{42}^3y, \zeta_{42}^{18}t) \\ \hline
36 & y^2=x^3-t^5(t^6+1)  & (x,y,t)\mapsto (\zeta_{36}^2x,\zeta_{36}^3y,\zeta_{36}^{30}t) \\ \hline
28 & y^2=x^3+x+t^7       & (x,y,t)\mapsto (-x, \zeta_{28}^7y,\zeta_{28}^2t) \\ \hline
12 & y^2=x^3+t^5(t^2+1)  & (x,y,t)\mapsto (\zeta_{12}^2x,\zeta_{12}^3y,-t) \\
\hline
\end{array}
\] 
\caption{Elliptic fibrations of unimodular K3 surfaces}
\end{table}
\medskip

The corresponding theorem in the non-unimodular case reads as follows:

\begin{Theorem}[Non-unimodular cases]
{\sl Let $X$ be an algebraic K3 surface over $\C$.

Assume that $T_X$ is non-unimodular. Let $\Omega=\{3, 9, 27, 5, 25, 7, 11, 13, 17, 19\}$. %\{3^m\, (1\leq m\leq 3),\, 5^{\ell}\,(1\leq \ell\leq 2),\, 7, 11, 13, 17, 19\}$. 
Denote $k=|H_X|$. Then the following assertions hold.

(a) Suppose that $\mbox{rank}(T_X)=\phi(k)$. Then $k\in\Omega$.

(b) Conversely, for each $k\in\Omega$, there exists a unique algebraic K3 surface
$X$ with $|H_X|=k$ and $\mbox{rank}(T_X)=\phi(k)$. The non-unimodular lattices $S_X, T_X$ are given in the following table.} 
\end{Theorem}

{\scriptsize
\begin{table}[ht!]
\[
\begin{array}{c||c|c||c|c} \hline
k & S_X & \mbox{rank}(S_X) & T_X& \mbox{rank}(T_X) \\ \hline\hline
19 & U_2\oplus \begin{pmatrix} -2 & 1 \\ 1 & -10\end{pmatrix} & 4 & (-E_8)^2\oplus\begin{pmatrix} 2 & 1 \\ 1 & 10\end{pmatrix} & 18 \\ \hline
17 & U_2\oplus{\scriptscriptstyle\begin{pmatrix} -2 & 0 & 0 & 1\\0 & -2 & 1 & 1\\0 & 1 & -2 & 0\\1 & 1 & 0 & -4\end{pmatrix}} & 6 & U_2^2\oplus(-E_8)\oplus\scriptscriptstyle{\begin{pmatrix} -2 & 0 & 0 & 1\\0 & -2 & 1 & 1\\0 & 1 & -2 & 0\\1 & 1 & 0 & -4\end{pmatrix}} & 16 \\ \hline
13 & (-E_8) \oplus  \begin{pmatrix} -2 & 5 \\ 5 & -6\end{pmatrix}  & 10 & U_2\oplus (-E_8) \oplus  \begin{pmatrix} -2 & 5 \\ 5 & -6\end{pmatrix} & 12 \\ \hline
11 & U_2\oplus(-A_{10})  & 12 & (-E_8)\oplus\begin{pmatrix} 2 & 1 \\ 1 & 6\end{pmatrix} & 10 \\ \hline
7 & U_2\oplus(-E_8)\oplus(-A_6) & 16 & U_2^2\oplus \begin{pmatrix} -2 & 1 \\ 1 & -4\end{pmatrix} & 6 \\ \hline
25  & \begin{pmatrix} -2 & 3 \\ 3 & -2\end{pmatrix} & 2 & U_2\oplus(-E_8)^2\oplus  \begin{pmatrix} -2 & 3 \\ 3 & -2\end{pmatrix} & 20 \\ \hline
5  & (-E_8)^2 \oplus  \begin{pmatrix} -2 & 3 \\ 3 & -2\end{pmatrix} & 18 & U_2\oplus  \begin{pmatrix} -2 & 3 \\ 3 & -2\end{pmatrix}  & 4 \\ \hline          
27 & U_2 \oplus (-A_2) & 4 & U_2^2\oplus (-E_6)\oplus (-E_8) & 18 \\ \hline
9 & U_2\oplus (-E_6)\oplus (-E_8) & 16 & U_2^2 + (-A_2) & 6 \\ \hline
3 & U_2\oplus A_2\oplus (-E_8)^2 & 20  & \begin{pmatrix} 2 & 1 \\ 1 & 2\end{pmatrix} & 2 \\
\hline
\end{array}
\]
\caption{Non-unimodular non-symplectic K3 surfaces}
\end{table}
}

\medskip
%
%{\bf Question:} Is there any divisibilty result as in Theorem 1? O-Z do not state it, but maybe...
%\marginpar{Q}

%{\bf Noriko's comment:} No, Vorontsov did not give any divisibility result.

\begin{Remark}
{\rm
In the general non-unimodular case, there is a divisibility result similar to (a) in Theorem 1:
Let $k=|H_X|$ without the assumption $\mbox{rank}(T_X)=\phi(k)$. Then
$k\in\Omega \cup \{1,2,4,8,16\}$ by \cite[Cor.~6.2]{Ko}. The cases where $k=2^j$ are analyzed by the second author in \cite{Sch2}.
}
\end{Remark}

\medskip

In a weaker form, this classification was announced by
Vorontsov \cite{Vo}. Kond\=o gave a proof of part (a) in \cite{Ko}. The uniqueness of part (b) was established by Machida and Oguiso \cite{MO} for $k=25$
and by Oguiso and Zhang \cite{OZ} for the remaining cases. The lattices will be calculated in the next section.

\medskip

Explicit defining equations for these K3 surfaces were given by Kond\=o \cite{Ko}.
All  K3 surfaces but one are elliptic with section and thus defined by Weierstrass equations.
The exception is the K3 surface corresponding to the case $k=25$. This is only defined as a double sextic. The following table reproduces the defining affine equations from \cite{Ko} up to sign changes. For the $H_X$-action we employ the same convention as before.

\begin{table}[ht!]
\[
\begin{array}{c|c|c} \hline
k & X & g_k \\ \hline\hline
19 & y^2=x^3+t^7x-t & (x,y,t)\mapsto (\zeta_{19}^7x, \zeta_{19}y, \zeta_{19}^2 t) \\ \hline
17 & y^2=x^3+t^7x-t^2  & (x,y,t)\mapsto (\zeta_{17}^7x, \zeta_{17}^2y, \zeta_{17}^2t) \\ \hline
13 & y^2=x^3+t^5x-t  & (x,y,t)\mapsto (\zeta_{13}^5x, \zeta_{13}y, \zeta_{13}^2t) \\ \hline
11 & y^2=x^3+t^5x-t^2  & (x,y,t)\mapsto (\zeta_{11}^5x,\zeta_{11}^2y,\zeta_{11}^2t) \\ \hline
7 & y^2=x^3+t^3x-t^8 & (x,y,t)\mapsto (\zeta_{7}^3x, \zeta_{7}y,\zeta_{7}^2t) \\ \hline
5 & y^2=x^3+t^3x-t^7 & (x,y,t)\mapsto (\zeta_5^3x,\zeta_5^2y,\zeta_5^2t) \\ \hline
27& y^2=x^3-t(t^9+1) & (x,y,t)\mapsto (\zeta_{27}^2x,\zeta_{27}^3y,\zeta_{27}^6t) \\ \hline
9 & y^2=x^3-t^5(t^3+1) & (x,y,t)\mapsto (\zeta_9^2x,\zeta_9^3y,\zeta_9^3t)  \\ \hline
3 & y^2=x^3+t^5(t-1)^2 & (x,y,t)\mapsto (\zeta_3x,y,t) \\ \hline
25& y^2=u^5+uv^5-1 & (u,v,y)\mapsto(\zeta_{25}^{20}u,\zeta_{25}v,y)\\
\hline
\end{array}
\]
\caption{Equations of non-unimodular K3 surfaces}
\end{table}
\medskip

\begin{Remark}
{\rm
 Though the K3 surfaces in Theorem 1 and Theorem 2 are unique
up to isomorphism, there are several ways of defining these K3 surfaces. For instance,
for $k=66$, we may take a weighted K3 surface: 
$$x_0^2+x_1^3+x_2^{11}x_3+x_3^{12}=0\subset\P(6,4,1,1)$$
of degree $12$. %This is unimodular. 
Letting $x_2=1$, we obtain the affine piece 
$x_0^2+x_1^3+x_3+x_3^{12}=0$. This is birationally equivalent to the elliptic surface defined by
Kond\=o, $y^2=x^3+t(t^{11}-1)$. This was also pointed out to us by Y. Goto.
We will elaborate on two examples in more detail in Remarks 5, 6.}
\end{Remark}

Some of the K3 surfaces in Theorem 1 and 2 offer geometric interpretations of the symplectic group actions. 
We elaborate on three cases in the next two examples.

\begin{Example}
{\rm 
The surfaces for $k=3$ and $k=12$ naturally occur in a one-dimensional family
\begin{eqnarray}\label{eq:Inose}
X_\lambda:\;\;\; y^2 = x^3 + t^5 (t^2 +2 \lambda t+1).
\end{eqnarray}
The K3 surfaces $X_\lambda$ are double coverings of the Kummer surfaces for $E_0\times E$ where $E_0$ is the elliptic curve with $j(E_0)=0$ and $E$ varies with $\lambda$ (cf.~Inose \cite{Inose}). The double covering exhibits the Shioda-Inose structure on $X_\lambda$. 

On the one hand, $X_\lambda$ admits an automorphism of order $k=12$ if and only if the elliptic curve $E$ admits an automorphism of order four, i.e.~$j(E)=1728$. This is the case $\lambda=0$. On the other hand, every $X_\lambda$ admits an automorphism of order $k=3$, since the elliptic fibration is isotrivial:
\[
g_3: x\mapsto \zeta_3\, x.
\]
Here $\mbox{rank}(T_{X_\lambda})=\phi(k)=2$ if and only if $E_0$ and $E$ are isogenous. However, unless $E\cong E_0$, the elliptic fibration (\ref{eq:Inose}) has Mordell-Weil rank two with non-trivial $g_k$-action. Hence  the only case with $\rho(X_\lambda)=20$ and trivial $g_k$-action on $S_{X_\lambda}$ is $E\cong E_0$. This corresponds to $\lambda=\pm 1$ as in the table.

By similar arguments, we can rule out the Kummer surfaces themselves, although they naturally inherit an automorphism of order three from $E_0$: Since the Mordell-Weil rank is always positive, this automorphism acts non-trivially on the N\'eron-Severi lattice.}
\end{Example}

\begin{Example}
{\rm
A similar picture involving Shioda-Inose structures arises for $k=5$: By \cite{Kumar}, $X$ admits a Nikulin involution $\iota$
\[
x\mapsto x^3/t^8,\;\;\; y\mapsto x^3\,y/t^{12},\;\;\; t\mapsto -x/t^3.
\]
The minimal resolution of the quotient $X/\iota$ is $\Q$-isomorphic to the Kummer surface of the Jacobian of the following genus 2 curve:
\[
C:\;\;\; u^2 = v^6 + 4 v. 
\]
Here $C$ is equipped with an automorphism of order five, $v\mapsto \zeta_5 v,\;u\mapsto \zeta_5^3 u$.}
\end{Example}

\end{section}

\begin{section}{Algebraic and transcendental cycles}
\label{s:cycles}

In this section, we compute the lattices $S_X, T_X$ of algebraic and transcendental cycles on the K3 surfaces $X$ from Theorem 2. We use the theory of the discriminant form as developed by Nikulin in \cite{N} and Mordell-Weil lattices after Shioda \cite{Sh-MW}.

\medskip

We start by computing the N\'eron-Severi lattices for the K3 surfaces in Theorem 2. For all but $k=25$, we will use the elliptic fibration given by Kond\=o \cite{Ko}.

\medskip

For $k=3,9,27$, there is no section. By the formula of Shioda and Tate, $S_X$ agrees with the trivial lattice which is generated by the zero section $O$, a general fiber $F$ and fiber components not meeting the zero section. After identifying the singular fibers with the corresponding Dynkin diagrams, we deduce the claimed shape.

\medskip

In the other cases, there is a non-torsion section $P$ that we give in the next table. Then the lattice $S_X$ is encoded in the intersection number $(P.O)$ and in the fiber components that $P$ meets. After Shioda \cite{Sh-MW}, this information can be expressed through the height of $P$: 

\[
h(P) = 4 + 2\,(P.O) - \sum_v \text{corr}_v(P).
\]

Here the sum runs over all reducible fibers and the correction terms behave as follows: On the one hand,

\[
\text{corr}_v(P) = 0\;\; \Leftrightarrow\;\; \text{$P$ meets the identity component of the fiber $F_v$}.
\]

\medskip

Otherwise we identify the singular fibers with the corresponding Dynkin diagrams. For $A_{n-1}$ (i.e.~$III, IV, I_n$ in Kodaira's notation), we number the components $\Theta_i$ cyclically such that $O$ meets $\Theta_0$. The entry in the following table for $A_{n-1}$ lists the correction term if $P$ meets $\Theta_j\, (j>0)$:

\medskip
\begin{table}[ht!]
$$
\begin{array}{c||c|c|c}
\hline
\text{fiber type} & E_6 & E_7 & A_{n-1}\\
\hline
\text{corr}_v(P) & 4/3 & 3/2 & j(n-j)/n\\
\hline
\end{array}
$$
\end{table}

\medskip

In \cite{Sh-MW}, the height was introduced to endow the Mordell-Weil group modulo torsion with the structure of a positive definite lattice (though not integral in general). For the K3 surfaces in consideration, the Mordell-Weil rank can only be one, since the ranks of $S_X$ and $T_X$ add up to 22. Moreover there cannot be torsion in the Mordell-Weil group. Hence the discriminant of $S_X$ is given by the following formula:

\[
\mbox{disc}(S_X) = h(P) \, \prod_v \mbox{disc}(F_v).
\]

\medskip

The following table lists the reducible singular fibers, the MW-generator and its height plus the resulting discriminant for the elliptic K3 surfaces in Theorem 2. Some sections involve the imaginary unit $i=\sqrt{-1}$.

\medskip

\begin{table}[ht!]
$$
\begin{array}{c||c|c|c|c}
\hline
k & \mbox{reducible fibers} & P & h(P) & \mbox{disc}(S_X)\\
\hline \hline
19 & III & (1/t^6, 1/t^9) & 19/2 & -19\\
\hline
17 & III, IV & (0,i\,t) & 17/6 & -17\\
\hline
13 & III^* & (1/t^4, 1/t^6) & 13/2 & -13\\
\hline
11 & IV, III^* & (0,i\,t) & 11/6 & -11\\
\hline
7 & IV^*, III^* & (0,i\,t^4) & 7/6 & -7\\
\hline
5 & III^*, II^* & (t^4, t^6) & 5/2 & -5\\
\hline
27 & IV & - & - & -3\\
\hline
9 & IV^*, II^* & - & - & -3\\
\hline
3 & IV, II^*, II^* & - & - & -3\\
\hline
\end{array}
$$
\caption{Reducible singular fibers and sections of non-unimodular elliptic K3 surfaces}
\label{Tab:5}
\end{table}

\medskip

For $k=17, 19$, the N\'eron-Severi lattice $S_X$ given in Theorem 2 is exactly the hyperbolic plane $U_2$ generated by $O$ and $F$, plus its orthogonal complement. The orthogonal projection $\pi$ takes $P$ to the divisor class of
\medskip
\[
\pi(P)=P-O-((P.O)+2) F \;\;\;\text{of self intersection}\;\;\;\pi(P)^2=-4-2 (P.O).
\]

\medskip

For the other surfaces, we shall use the discriminant form to determine the abstract shape of the lattices. This approach will also suffice to find the lattices $T_X$ of transcendental cycles for all surfaces from Theorem 2.

\medskip

Given an even integral non-degenerate lattice $L$, we denote its dual by $L^\vee$. In \cite{N}, Nikulin introduced a quadratic form  on the quotient $L^\vee/L$ which he called discriminant form:

\begin{eqnarray*}
q_L:\;\; L^\vee/L & \to & \Q \mod 2\Z\\
x & \mapsto & x^2 
\end{eqnarray*}

\medskip

\begin{Theorem}[Nikulin {\cite[Cor.~1.9.4]{N}}]
The genus of an even integral non-degenerate lattice is determined by its signature and discriminant form.
\end{Theorem}

\medskip

Consider the above cases of Mordell-Weil rank one ($k=5,7,11,13,17,19$). Then $L^\vee/L$ is cyclic of order $k$. The discriminant form maps the canonical generator to $-1/h(P)$. This is abbreviated by the notation

\[
q_{S_X} = \Z/k\Z\left(-\frac{1}{h(P)}\right).
\]

\medskip

In each case, it is immediate that the claimed lattice in Theorem 2 has exactly the same discriminant form. By \cite[\S 15, Cor.~22]{CS}, the discriminant has small enough absolute value so that there is only one class per genus.  Hence the lattices are isomorphic by Theorem 3.

\medskip

In the non-elliptic case $k=25$, we proceed as follows: Since $T_X$ has rank 20, we know that $\rho(X)=2$. We find generators of $S_X$ on the double sextic model from Theorem 2. Consider the line 
\[
\ell=\{u=0\}\subset\P^2.
\]
This line meets the branch curve of the double cover with multiplicity six at one point. The pull-back $\pi^*\ell$ to the double cover $X$ splits into two rational curves
\[
\ell_\pm=\{u=0, y=\pm \sqrt{-1}\}\subset X.
\]
These lines intersect with multiplicity three at the preimage of the above point. Since $(\pi^*\ell)^2=2$, we deduce $\ell_\pm^2=-2$. This gives the intersection matrix from Theorem 2. As its determinant $-5$ is squarefree and $\rho(X)=2$, we have $S_X=\langle \ell_\pm\rangle$.

\medskip

\begin{Example}
{\rm
For $k=7$, there is a geometric way to see the abstract shape of the N\'eron-Severi lattice $S_X$: We find an elliptic fibration on $X$ with trivial lattice $S_X$. In terms of the Weierstrass equation above, we blow up three times at $(0,0,0)$. One chart then is given by
\[
x=t^3 x',\;\;\; y=t^3 y'.
\]
After dividing by the common multiple, the resulting equation is
\[
y'^2 = x'^3 t^3 + t^2 + x'.
\]
We want to use $x'=x/t^3$ as the coordinate of the base curve $\P^1$. A change of variables gives rise to the Weierstrass form
\[
y'^2 = t^3 + t^2 + x'^7.
\]
This has the claimed trivial lattice $U_2\oplus(-A_6)\oplus(-E_8)$ with trivial action of the non-symplectic automorphism $x'\mapsto \zeta_7\, x'$.}
\end{Example}

\begin{Example}
{\rm
A similar approach works for $k=11$. By \cite[Lem.~2.1, 2.2]{Ko}, the surface $X$ from Theorem 2 admits an elliptic fibration with section and a singular fiber of type $I_{11}$.
We claim that this fibration is given as follows
\[
Z:\;\;\; y^2 = x^3 + x^2 + t^{11}.
\]
By definition, $Z$
has the claimed trivial lattice $U_2\oplus(-A_{10})$ with trivial operation of the non-symplectic automorphism
\[
\varphi: \;\; t\mapsto \zeta_{11} \, t.
\]
However, we did not find an explicit transformation between the above Weierstrass form and the one from Theorem 2. Instead we can use the uniqueness of Theorem 2 to prove that the K3 surfaces are isomorphic. 

To see this, we only need that $\phi(11)=10\, |\, \mbox{rank}(T_Z)$. Since $T_Z$ has rank $22-\rho(Z)\leq 10$, we deduce equality and $\rho(Z)=12$. In particular, $\varphi$ operates trivially on $S_Z$. By uniqueness,  $Z$ is isomorphic over $\C$ to the K3 surface from Theorem 2.}

%For this we only need that $MW(Z)=\{O\}$, since then $S_Z=U_2\oplus(-A_{10})$. We shall now assume that there is a non-trivial section $P\in MW(Z)$ and establish a contradiction.

%\medskip

%If $P$ was $\varphi$-invariant, then it would descend to the rational elliptic surface
%\[
%W:\;\;\; y^2 = x^3 + x^2 + t.
%\]
%However, $W$ has trivial Mordell-Weil group since $\rho(W)=10$ and there is a fibre of type $II^*$ at $\infty$. Hence $MW(Z)^\varphi=\{O\}$. Then the $\Z[\varphi]$-module structure on $MW(Z)$ implies that
%\[
%\mbox{rank}(MW(Z)) \equiv 0 \mod 10.
%\]
%Hence $\rho(Z)=12+\mbox{rank}(MW(Z))\equiv 2 \mod 10$. Since $\rho(Z)\leq h^{1,1}(Z)=20$, we deduce $\rho(Z)=12$ and $MW(Z)=\{O\}$. Hence $\varphi$ acts trivially on $S_Z$. The uniqueness in Theorem 2 implies the claim.
\end{Example}

It remains to prove the lattices $T_X$ of transcendental cycles in Theorem 2. This only requires the following result:

\begin{Theorem}[Nikulin {\cite[Prop.~1.6.1]{N}}]
Let $N$ be an even integral unimodular lattice. Let $L$ be a primitive non-degenerate sublattice and $M=L^\bot$. Then
\[
q_L = -q_M.
\]
\end{Theorem}

On a K3 surface $X$, the N\'eron-Severi lattice $S_X$ always embeds primitively into $H^2(X,\Z)$. Hence the proof of the transcendental lattices in Theorem 2 is an easy application of Theorem 4. We use that there is one class per genus in each case. For $k\neq 3$, $T_X$ is indefinite, so we deduce the claim from \cite[\S 15, Cor.~22]{CS}. For $k=3$, $T_X$ is positive definite, so the claim follows from class group theory.

\end{section}

\begin{section}{Delsarte surfaces}
The K3 surfaces listed in Theorem 1 and Theorem 2 are all Delsarte surfaces except for $k=3$. Indeed,
each surface is defined by a sum of four monomials in the affine $3$-space. %, or in a weighted projective $3$-space. 
Hence such a surface is covered by a Fermat surface of a suitable degree.
The main result of this section is to determine the corresponding Fermat surface
and a covering map for each of these K3 surfaces.

From now on, we let 
\begin{eqnarray}\label{eq:F_m}
\mathcal{F}_m\,:\, x_0^m+x_1^m+x_2^m+x_3^m=0\subset\P^3
\end{eqnarray}
denote the Fermat surface
of degree $m\geq 4$.  
\medskip

\begin{Theorem} \label{Thm:3}
{\sl Let $X$ be one of the K3 surfaces listed in Theorem 1 and Theorem 2. Let $k=|H_X|$.
If $k\neq 3$, then $X$ is covered by the Fermat surface of the following degree $m$:
$$
\begin{cases}
m=k, & \text{{\sl if  $X$ is unimodular,}}\\
m=2k, & \text{{\sl if  $X$ is non-unimodular.}}
\end{cases}
$$}
\end{Theorem}

For each K3 surface, we will prove Theorem \ref{Thm:3} by explicitly giving the covering map on the affine model $\{x_0\neq 0\}$ of $\mathcal{F}_m$ with coordinates $U=\frac{x_1}{x_0}, V=\frac{x_2}{x_0}, W=\frac{x_3}{x_0}$:
\[
\mathcal{F}_m:\;\;\; U^m + V^m + W^m + 1 = 0.
\]
A method to compute these covering maps goes back to Shioda in \cite{Sh}.

\bigskip

\noindent{\bf (I) Unimodular cases}
%\smallskip
\medskip

\noindent$\boxed{k=66}$ The elliptic surface is given by the Weierstrass equation
$$y^2=x^3-t(t^{11}+1)=x^3-t^{12}-t.$$
It is more convenient to work in the chart at $\infty$ where $t\neq 0$. The local parameter is $s=1/t$. Then the equation is transformed to
$$y^2=x^3-1/s^{12}-s^{11}/s^{12}.$$
Introduce new coordinates $\xi:=s^4\,x,\, \eta:=s^6\,y$. The
resulting equation is
$$\eta^2=\xi^3-1-s^{11}.$$
The covering map from the affine Fermat surface of degree $66$ is as follows:
%$$U^{66}=V^{66}-1+W^{66}\mapsto \eta^2=\xi^3+1-s^{11}$$
$$\eta\mapsto U^{33},\; \xi\mapsto -V^{22},\; s\mapsto W^{6}.$$

\medskip

\noindent$\boxed{k=44}$
The covering map from the affine Fermat surface of degree $44$ is as follows:
$$\{U^{44}+V^{44}+W^{44}+1=0\}\to \{y^2=x^3+x+t^{11}\}$$
$$y\mapsto U^{22}\,V^{11},\, x\mapsto -V^{22},\, t\mapsto -W^{4}\,V^2.$$

\medskip

\noindent$\boxed{k=42}$  The argument is very similar to the
case $k=66$. The same transformations lead to
$$\eta^2=\xi^3-1-s^7.$$
The covering map from the affine Fermat surface of degree $42$ is as follows:
%$$\{U^{42}+V^{42}+1-W^{42}\}\to \{\eta^2=\xi^3+1-s^7\}$$
$$\eta\mapsto U^{21},\; \xi\mapsto -V^{14},\; s\mapsto W^6.$$
\medskip

\noindent$\boxed{k=36}$
The covering map from the affine Fermat surface of degree $36$ is as follows:
%$$\{U^{36}=V^{36}-W^{36}+1=0\}\to y^2=x^3-t^5(t^6-1)$$
$$y\mapsto U^{18}\,W^{15},\; x\mapsto -V^{12}\,W^{10},\; t\mapsto W^6.$$

\medskip

\noindent$\boxed{k=28}$
The covering map from the affine Fermat surface of degree $28$ is as follows:
%$$U^{28}=V^{28}+W^{28}+1\mapsto y^2=x^3+x+t^7$$
$$y\mapsto U^{14}\,V^7,\; x\mapsto -V^{14},\; t\mapsto -W^4\,V^2.$$

\medskip

\noindent$\boxed{k=12}$ 
The covering map from the affine Fermat surface of degree $12$ is as follows:
%$$U^{12}=V^{12}+W^{12}+1\mapsto y^2=x^3+t^5(t^2+1)$$
$$y\mapsto U^{6}\,W^{15},\; x\mapsto -V^{4}\,W^{10},\; t\mapsto -W^6.$$

\medskip
%\vskip 1cm

%\pagebreak

\noindent{\bf (II) Non-unimodular cases}
%\smallskip
\medskip

\noindent$\boxed{k=19}$
The covering map from the affine Fermat surface of degree $38$ is as follows:
%$$U^{38}=V^{38}+W^{38}-1\mapsto y^2=x^3+t^7x-t$$
$$y\mapsto U^{19}\,W^3/V,\; x\mapsto -V^{12}\,W^2,\; t\mapsto W^6/V^2.$$

\medskip

\noindent$\boxed{k=17}$
The covering map from the affine Fermat surface of degree $34$ is as follows:
%$$U^{34}=V^{34}+W^{34}+1\mapsto y^2=x^3+t^7x+t^2$$
$$y\mapsto U^{17}\,W^6/V^2,\; x\mapsto -V^{10}\,W^4,\; t\mapsto W^6/V^2.$$

\medskip

\noindent$\boxed{k=13}$
The covering map from the affine Fermat surface of degree $26$ is as follows:
%$$U^{26}=V^{26}+W^{26}-1\mapsto y^2=x^3+t^5x-t$$
$$y\mapsto U^{13}\,W^3/V,\; x\mapsto -V^{8}\,W^2,\; t\mapsto W^6/V^2.$$

\medskip

\noindent$\boxed{k=11}$
The covering map from the affine Fermat surface of degree $22$ is as follows:
%$$U^{22}=V^{22}+W^{22}+1\mapsto y^2=x^3+t^5x+t^2$$
$$y\mapsto U^{11}\,W^6/V^2,\; x\mapsto -V^{6}\,W^4,\; t\mapsto W^6/V^2.$$

\medskip

\noindent$\boxed{k=7}$
The covering map from the affine Fermat surface of degree $14$ is as follows:
%$$U^{14}=V^{14}+W^{14}+1\mapsto y^2=x^3+t^3x+t^8$$
$$y\mapsto U^{7}\,V^8/W^{24},\; x\mapsto -V^{10}/W^{16},\; t\mapsto V^2/W^6.$$

\medskip

\noindent$\boxed{k=5}$
The covering map from the affine Fermat surface of degree $10$ is as follows:
%$$U^{10}=V^{10}+W^{10}-1\mapsto y^2=x^3+t^3x-t^7$$
$$y\mapsto U^{5}\,V^7/W^{21},\; x\mapsto -V^{8}/W^{14},\; t\mapsto V^2/W^6.$$

\medskip

\noindent$\boxed{k=27}$ 
The covering map from the affine Fermat surface of degree $54$ is as follows:
%$$U^{54}=V^{54}+W^{54}-1\mapsto y^2=x^3+t(t^9-1)$$
$$y\mapsto U^{27}\,W^{3},\; x\mapsto -V^{18}\,W^{2},\; t\mapsto W^6.$$

\medskip

\noindent$\boxed{k=9}$ 
The covering map from the affine Fermat surface of degree $18$ is as follows:
%$$U^{18}=V^{18}+W^{18}-1\mapsto y^2=x^3+t^5(t^3-1)$$
$$y\mapsto U^{9}\,W^{15},\; x\mapsto -V^{6}\,W^{10},\; t\mapsto W^6.$$

%\medskip

%\noindent$\boxed{k=3}$  
%The covering map from the affine Fermat surface of degree $12$ is as follows:
%$$U^{12}=V^{12}+W^{12}-1\mapsto y^2=x^3+(t-1)^2$$
%$$\eta\mapsto U^{15},\; \xi\mapsto V^{10},\; s\mapsto W^3.$$

\medskip

\noindent$\boxed{k=25}$
Affinely, the covering map from the Fermat surface of degree $50$ is as follows:
%$$U^{50}=V^{50}+W^{50}+1\mapsto y^2=1+u^5+uv^5$$
$$y\mapsto U^{25},\; u\mapsto -V^{10},\; v\mapsto W^{10}/V^2.$$

\end{section}

\begin{section}{Motivic decomposition}
\label{s:motif}

Shioda studied Delsarte surfaces in \cite{Sh}. He computed their Picard numbers through the covering Fermat surfaces. Here we recall his argument and apply it to the K3 surfaces in Theorem 1 and 2. Note that we will actually be most interested in the transcendental lattices.

\medskip

Let $\bmu_m$ be the
group of $m$-th roots of unity and $\Delta$ the image of the diagonal
inclusion $\bmu_m\hookrightarrow\bmu_m^4$. Then the quotient group $M=\bmu_m^{4}/\Delta$ operates by multiplication of coordinates on the Fermat surface $\mathcal{F}_m$.

\medskip

This group operation induces a decomposition of the second cohomology $H^2(\mathcal{F}_m)$ 
into one-dimensional eigenspaces with character. More generally true for Fermat varieties 
of any dimension, this result is due to Weil \cite{We1}. The decomposition is best described in terms 
of the character group of $M$:
$$\mathfrak{A}_m:=
\{\alpha=(a_0,a_1,a_2,a_3)\in(\Z/m\Z)^4\,|\, a_i\not\equiv 0\pmod m,\,\sum_{i=0}^3 a_i\equiv 0\pmod m\,\}.
$$
It is well-known that $H^2(\mathcal{F}_m)$ decomposes into the subspace $V_0$ of the hyperplane class $H$ and one-dimensional subspaces $V(\alpha)$ for each $\alpha\in\mathfrak{A}_m$:
\begin{eqnarray}\label{eq:decomp}
H^2(\mathcal{F}_m) = V_0 \oplus \bigoplus_{\alpha\in\mathfrak{A}_m} V(\alpha).
\end{eqnarray}

%\begin{Remark}
%\label{Rem7}
%{\rm 1. Each character space in \eqref{eq:decomp} is one dimensional. Hence the $\Gal(\Qbar/\Q(\zeta_k))$
%representation $H^2(\mathcal{F}_m \times \Q(\zeta_k))$ is a sum of one dimensional representations. In fact
%Weil proved that each of these representations is a Gr\"ossencharactere associated to the Jacobi sum
%determined by the character $\alpha$. In particular, the $V(\alpha)$ are pairwise nonisomorphic as 
%$\Gal(\Qbar/\Q(\zeta_k))$ 
%representations.\\
%2. The action of  
%$t\in\Gal(\Q(\zeta_k))/\Q\simeq (\Z/k\Z)^\ast$ 
%sends $V(\alpha)$ to $V(t\alpha)$.\\
%3. The character $\alpha$ contributes to the Hodge group $H^{|\alpha|-1,3-|\alpha|}(\mathcal{F}_m)$.
%}
%\end{Remark}
We now describe a criterion whether $V(\alpha)$ is algebraic: If we choose representatives $0<a_i<m$, then this gives a well-defined map
\[
|\alpha|=\frac 1m \sum_{i=0}^3 a_i.
\]
Let $u\in(\Z/m\Z)^*$ operate on $\alpha\in\mathfrak{A}_m$ coordinatewise by multiplication. Then we define a subset $\mathfrak{B}_m$ of $\mathfrak{A}_m$ as follows:
\[
\mathfrak{B}_m = \{\alpha\in \mathfrak{A}_m; |u\cdot\alpha|=2 \;\; \forall u\in(\Z/m\Z)^*\}.
\]
{\bf Criterion:} 
{\sl The eigenspace $V(\alpha)$ is algebraic if and only if $\alpha\in\mathfrak{B}_m$.}

\medskip

\begin{Remark}\label{Rem7}
{\rm 
The criterion is based on the following two facts:
\begin{enumerate}
\item The action of  
$u\in\Gal(\Q(\zeta_m)/\Q)\simeq (\Z/m\Z)^\ast$ 
sends $V(\alpha)$ to $V(u\cdot \alpha)$.
\item\label{rem-2}
The eigenspace $V(\alpha)$ with character $\alpha\in\mathfrak{A}_m$ contributes to the Hodge group $H^{|\alpha|-1,3-|\alpha|}(\mathcal{F}_m)$.
\end{enumerate}}
\end{Remark}

\medskip

Using the above criterion, one can easily determine the  transcendental and algebraic part of $H^2(\mathcal{F}_m)$ and in particular the Picard number. We shall now see how this carries over to the K3 surfaces in Theorem 1 and 2. %Most of these arguments are based on Shioda's work \cite{Sh}.

\begin{Lemma}\label{Lem:G}
{\sl Each K3 surface from Theorem 1 and 2 except for $k=3$ is birationally
equivalent to a quotient of the covering Fermat surface in Theorem \ref{Thm:3}.}
\end{Lemma}

{\sl Proof:}
In the proof of Theorem 3, we have exhibited the explicit covering map $\pi$ between the K3 surface $X$ and  the Fermat surface $\mathcal{F}_m$. We define the following subgroup $G$ of the automorphism group $M=\bmu_m^{4}/\Delta$ (viewed as $\bmu_m^3$ operating on the affine coordinates of $\mathcal{F}_m$):
\[
G=\{g\in M; \pi=\pi\circ g\}.
\]
It is immediate that $X$ is birationally equivalent to the quotient $\mathcal{F}_m/G$. \qed

\begin{Example}
{\sl In some cases, the precise shape of $G$ is visible directly. For instance, if $k=66$, then
\[
G=\bmu_{33}\times\bmu_{22}\times\bmu_6.
\]
In other cases such as $k=44$, some more work is required to determine $G$ abstractly, but we will not go into the details here.}
\end{Example}

Let $X$ be one of the K3 surfaces from Theorem 1 or 2 with $k\neq 3$ and $\mathcal{F}_m$ the covering Fermat surface. The key property for the analysis of $H^2(X)$ lies in the Lefschetz number $\lambda(X)$ which gives the dimension of the transcendental subspace of $H^2(X)$:
\[
\lambda(X)=b_2(X)-\rho(X).
\]
Namely, $\lambda(X)$ is a birational invariant of algebraic surfaces. Hence
\[
\lambda(X)=\lambda(\mathcal{F}_m/G).
\]
Thus we can identify the transcendental part of $H^2(X)$ with the transcendental part of $H^2(\mathcal{F}_m)$ that is invariant under $G$. Here we shall further use the decomposition (\ref{eq:decomp}) of  $H^2(\mathcal{F}_m)$:

\medskip

{\bf Criterion:} 
{\sl Let $\alpha=(a_0,a_1,a_2,a_3)\in\mathfrak{A}_m$. Then $V(\alpha)$ is $G$-invariant if and only if 
\[
\prod_{i=1}^3 \zeta_i^{a_i} =1 \;\;\;\forall\, g=(\zeta_1, \zeta_2, \zeta_3)\in G.
\]}

Let $\mathfrak{S}_G=\{\alpha\in\mathfrak{A}_m; V(\alpha) \text{ is $G$-invariant}\}$. Then
\[
H^2(\mathcal{F}_m)^G = V_0 \oplus \bigoplus_{\alpha\in\mathfrak{S}_G} V(\alpha).
\]
In consequence, the transcendental part $T(X)$ of $H^2(X)$ can be identified with
\begin{eqnarray}\label{eq:mot_dec}
T(X) = \bigoplus_{\alpha\in\mathfrak{S}_G\setminus\mathfrak{B}_m} V(\alpha).
\end{eqnarray}

In the following tables, we list the character sets $\mathfrak{S}_G\setminus\mathfrak{B}_m$ for all K3 surfaces from Theorem 1 and 2 except for $k=3$. We use the shorthand $\alpha=[a_1, a_2, a_3]$ for $\alpha=(a_0, a_1, a_2, a_3)$, corresponding to the affine chart $\{x_0\neq 1\}$, since this determines $a_0$ uniquely.

\begin{table}[ht!]
{\small
\[
\begin{array}{c|c} \hline \hline
k & \mathfrak{S}_G\setminus\mathfrak{B}_m\\ \hline\hline
66 & [6, 33, 22], [6, 33, 44], [12, 33, 22], [12, 33, 44], [18, 33, 22], [18, 33, 44], [24, 33, 22],\\
&  [24, 33, 44], [30, 33, 22], [30, 33, 44], [36, 33, 22], [36, 33, 44], [42, 33, 22], [42, 33, 44],\\
&  [48, 33, 22], [48, 33, 44], [54, 33, 22], [54, 33, 44], [60, 33, 22], [60, 33, 44]\\ 
\hline
44 & [1, 22, 24], [3, 22, 28], [5, 22, 32], [7, 22, 36], [9, 22, 40], [13, 22, 4], [15, 22, 8],\\
&  [17, 22, 12], [19, 22, 16], [21, 22, 20], [23, 22, 24], [25, 22, 28], [27, 22, 32], [29, 22, 36],\\
&  [31, 22, 40], [35, 22, 4], [37, 22, 8], [39, 22, 12], [41, 22, 16], [43, 22, 20]\\
\hline
42 & [6, 21, 14], [6, 21, 28], [12, 21, 14], [12, 21, 28], [18, 21, 14], [18, 21, 28],\\
& [24, 21, 14], [24, 21, 28], [30, 21, 14], [30, 21, 28], [36, 21, 14], [36, 21, 28]\\
\hline
36 & [1, 18, 12], [5, 18, 24], [7, 18, 12], [11, 18, 24], [13, 18, 12], [17, 18, 24],\\
&  [19, 18, 12], [23, 18, 24], [25, 18, 12], [29, 18, 24], [31, 18, 12], [35, 18, 24]\\
\hline
28 & [1, 14, 16], [3, 14, 20], [5, 14, 24], [9, 14, 4], [11, 14, 8], [13, 14, 12],\\
&  [15, 14, 16], [17, 14, 20], [19, 14, 24], [23, 14, 4], [25, 14, 8], [27, 14, 12]\\
\hline
12 & [1, 6, 4], [5, 6, 8], [7, 6, 4], [11, 6, 8]\\
\hline\hline
\end{array}
\]}
\caption{Motivic decomposition for unimodular K3 surfaces}
\end{table}

\begin{table}[ht!]
{\small
\[
\begin{array}{c|c} \hline \hline
k & \mathfrak{S}_G\setminus\mathfrak{B}_m\\ \hline\hline
19 & [19, 1, 35], [19, 3, 29], [19, 5, 23], [19, 7, 17], [19, 9, 11], [19, 11, 5],\\
&  [19, 13, 37], [19, 15, 31], [19, 17, 25], [19, 21, 13], [19, 23, 7], [19, 25, 1],\\
&  [19, 27, 33], [19, 29, 27], [19, 31, 21], [19, 33, 15], [19, 35, 9], [19, 37, 3]\\
\hline
17 & [17, 2, 28], [17, 4, 22], [17, 6, 16], [17, 8, 10], [17, 10, 4], [17, 12, 32],\\
&  [17, 14, 26], [17, 16, 20], [17, 18, 14], [17, 20, 8], [17, 22, 2],\\
&  [17, 24, 30], [17, 26, 24], [17, 28, 18], [17, 30, 12], [17, 32, 6]\\
\hline
13 & [13, 1, 23], [13, 3, 17], [13, 5, 11], [13, 7, 5], [13, 9, 25], [13, 11, 19],\\
&  [13, 15, 7], [13, 17, 1], [13, 19, 21], [13, 21, 15], [13, 23, 9], [13, 25, 3]\\
\hline
11 & [11, 2, 16], [11, 4, 10], [11, 6, 4], [11, 8, 20], [11, 10, 14],\\
&  [11, 12, 8], [11, 14, 2], [11, 16, 18], [11, 18, 12], [11, 20, 6]\\
\hline
7 & [7, 2, 8], [7, 4, 2], [7, 6, 10], [7, 8, 4], [7, 10, 12], [7, 12, 6]\\
\hline
25 &  [25, 2, 40], [25, 4, 30], [25, 6, 20], [25, 8, 10], [25, 12, 40], [25, 14, 30], [25, 16, 20],\\
&  [25, 18, 10], [25, 22, 40], [25, 24, 30], [25, 26, 20], [25, 28, 10], [25, 32, 40], [25, 34, 30],\\
&  [25, 36, 20], [25, 38, 10], [25, 42, 40], [25, 44, 30], [25, 46, 20], [25, 48, 10]\\
\hline
5 & [5, 1, 7], [5, 3, 1], [5, 7, 9], [5, 9, 3]\\
\hline
27 & [1, 36, 27], [5, 18, 27], [7, 36, 27], [11, 18, 27], [13, 36, 27], [17, 18, 27],\\
&  [19, 36, 27], [23, 18, 27], [25, 36, 27], [29, 18, 27], [31, 36, 27], [35, 18, 27],\\
&  [37, 36, 27], [41, 18, 27], [43, 36, 27], [47, 18, 27], [49, 36, 27], [53, 18, 27]\\
\hline
9 & [1, 6, 9], [5, 12, 9], [7, 6, 9], [11, 12, 9], [13, 6, 9], [17, 12, 9]\\
\hline\hline
\end{array}
\]}
\caption{Motivic decomposition for non-unimodular K3 surfaces}
\end{table}

\begin{Remark}
\label{Rem8}
{\rm 1. For each $k$, it is easily checked that 
the set $\mathfrak{S}_G\setminus\mathfrak{B}_m$
constitutes a single $(\Z/m\Z)^{\times}$-orbit.\\
2. We also see by Remark \ref{Rem7}.2 that the first character in each 
entry is the unique one of Hodge type $(0,2)$.}
\end{Remark} 
\end{section}

\begin{section}{Modularity}

In this section, we will prove the modularity of all K3 surfaces from Theorem 1 and 2. We compute their $\zeta$-functions over finite fields explicitly in terms of Jacobi sums.

\medskip

For $k=3$, the K3 surface is singular ($\rho=20$). Hence its modularity follows from a result by the first author \cite{Li} (cf.~Rem.~\ref{Rem:k=3}). For all other K3 surfaces from Theorem 1 and 2, we will use the covering Fermat surface from section 3 and the motivic decomposition from section 4 to prove modularity in Lemma \ref{Lem:R_t}.

\medskip

It goes back to Weil \cite{We} that the Fermat surface $\mathcal{F}_m$ is of CM type in the following sense: Over the cyclotomic field $\Q(\zeta_m)$, the Galois representation of $H^2(\mathcal{F}_m)$ splits into one-dimensional subrepresentations corresponding to the eigenspaces $V_0, V(\alpha)\; (\alpha\in\mathfrak{A}_m)$ from section 4. These subrepresentations are associated to Hecke Gr\"ossencharacters and can be described in terms of Jacobi sums. Here the Jacobi sums are determined by the characters $\alpha$. 
In particular, the $V(\alpha)$ are pairwise non-isomorphic as 
$\Gal(\Qbar/\Q(\zeta_m))$ 
representations.

%In fact, they are in correspondence with the one-dimensional eigenspaces $V_0, V(\alpha)\; (\alpha\in\mathfrak{A}_m)$ from section 4. Thus these eigenspaces give the motivic decomposition of $H^2(\mathcal{F}_m)$.

%{\rm 1. Each character space in \eqref{eq:decomp} is one dimensional. Hence the $\Gal(\Qbar/\Q(\zeta_k))$
%representation $H^2(\mathcal{F}_m \times \Q(\zeta_k))$ is a sum of one dimensional representations. In fact
%Weil proved that each of these representations is a Gr\"ossencharactere associated to the Jacobi sum
%determined by the character $\alpha$. In particular, the $V(\alpha)$ are pairwise nonisomorphic as 
%$\Gal(\Qbar/\Q(\zeta_k))$ 
%representations.\\

\medskip

By Lemma \ref{Lem:G} (and by \cite{Li} for $k=3$), we immediately obtain the

\begin{Corollary}\label{Cor:CM}
{\sl The K3 surfaces in Theorem 1 and 2 are also of CM type.} 
\end{Corollary}

In fact, we can even determine the $\zeta$-function of all K3 surfaces in Theorem 1 and 2. For this, we recall Weil's theorem \cite{We} for the covering Fermat surface. To state the result, we recall the basic setup involving Jacobi sums.

\medskip

We fix the degree $m$ and a prime $p$. Let $q=p^r$ such that
\[
q\equiv 1 \mod m.
\]
On the field $\F_q$ of $q$ elements, we fix a character 
\[
\chi: \F_q^* \to \C^*
\]
of order exactly $m$. For any $\alpha=(a_0,a_1,a_2,a_3)\in\mathfrak{A}_m$, we then define the Jacobi sum
\[
j(\alpha) = \sum_\text{\small $\begin{matrix} v_1, v_2, v_3\in\F_q^*\\ v_1+v_2+v_3=-1\end{matrix}$} \chi(v_1)^{a_1} \chi(v_2)^{a_2}\chi(v_3)^{a_3}.
\]

\begin{Theorem}[Weil]
\label{Thm:Weil}
{\sl In the above notation, the Fermat surface $\mathcal{F}_m$ over $\F_q$ has the following $\zeta$-function:
\[
\zeta(\mathcal{F}_m/\F_q, T) = \dfrac 1{(1-T)\, P(T)\, (1-q^2\, T)}
\]
where
\[
P(T) = (1-q\,T)  \prod_{\alpha\in\mathfrak{A}_m} (1-j(\alpha)\,T).
\]}
\end{Theorem}

We will now use Theorem \ref{Thm:Weil} to determine the $\zeta$-functions of all K3 surfaces from Theorem 1 and 2. Since every K3 surface $X$ has $b_1(X)=0$,
\begin{eqnarray}\label{eq:R(T)}
\zeta(X/\F_q, T) = \dfrac 1{(1-T)\, R(T)\, (1-q^2\, T)}
\end{eqnarray}
where $R(T)$ is the reciprocal characteristic polynomial of Frob$_q^*$ on $H^2(X)$. (Here we use without specifying an appropriate Weil cohomology, say $\ell$-adic \'etale cohomology for some prime $\ell\neq p$ on the base change $\bar X$ over the algebraic closure $\bar\F_q$.)

\medskip

Over $\Q$, the polynomial $R(T)$ factors according to algebraic and transcendental part of $H^2(X)$:
\[
R(T) = R_a(T)\,R_t(T).
\]
Thanks to Theorem \ref{Thm:Weil}, the factor $R_t(T)$ corresponding to the transcendental part is determined by Lemma \ref{Lem:G} and the motivic decomposition (\ref{eq:mot_dec}). \medskip

\begin{Lemma}\label{Lem:R_t}
{\sl Let $X$ be one of the K3 surfaces from Theorem 1 or 2 with $k\neq 3$. Fix the above setup over $\F_q$ with $q\equiv 1 \mod m$. Then
\[
R_t(T) = \prod_{\alpha\in\mathfrak{S}_G\setminus\mathfrak{B}_m} (1-j(\alpha)\,T).
\]}
\end{Lemma}

\medskip

\begin{Remark}\label{Rem:k=3}
{\rm 
The singular K3 surface from Theorem 2 with $k=3$ is modular by \cite{Li}. The affine model in section 2 is associated to the modular form of weight 3 and level 27 given in \cite[Table 1]{Sch1}. Over $\F_p\; (p\neq 2,3)$, the reciprocal characteristic polynomial $R_t(T)$ is given as follows:
\[
R_t(T) =
\begin{cases}
1-(\pi^2+\bar\pi^2) T + p^2 T^2, & \text{if $p=\pi \bar\pi$ in } \Z[3\,\zeta_3],\\
1 - p^2 T^2, & \text{if $p$ is inert in $\Q(\sqrt{-3})$}.
\end{cases}
\]
Because of this explicit description, we did not check how to express $R_t(T)$ in terms of Jacobi sums.}
\end{Remark}

\medskip

For the $\zeta$-functions (\ref{eq:R(T)}) of the K3 surfaces $X$ in Theorem 1 and 2, it remains to determine the reciprocal characteristic polynomial $R_a(T)$ of Frob$_q^*$ on the algebraic part of $H^2(X)$, i.e.~on $S_X$.

\medskip

\begin{Lemma}\label{Lem:NS}
{\sl Let $X$ a K3 surface in Theorem 1 or 2. Consider the affine model given in section 2 over some field $K$ of characteristic coprime to $2k$. Then $S_X$ is generated by algebraic cycles over $K(\sqrt{-1})$.}
\end{Lemma}

\medskip

{\sl Proof:}
We prove the claim for $K=\Q$. The lemma then follows by smooth base change.

\medskip

If $k\neq 25$, then we have already studied an elliptic fibration on $X$. It is easy to see that all fiber components are defined over $\Q(\sqrt{-1})$. Since the same holds for the sections by Table \ref{Tab:5}, the claim follows.

\medskip

For $k=25$, we saw that $S_X$ is generated by two lines which are conjugate in $\Q(\sqrt{-1})$. This proves the Lemma for $K=\Q$ and consequently for any field $K$ such that $X$ defines a (smooth) K3 surface over $K$. \qed

\medskip

Some non-unirational surfaces really require the extension by $\sqrt{-1}$. Explicitly, the local Euler factors take the following shape:

\begin{Corollary}\label{Cor:NS}
{\sl Let $X$ a K3 surface in Theorem 1 or 2. Consider the affine model from section 2 over some finite field $\F_q$. Let
\[
n_- = 
\begin{cases}
0, & q\equiv 1 \mod 4 \text{ or } k\neq 7,9,11,17,25,27,\\
1, & q\equiv 3 \mod 4 \text{ and } k=25, 27,\\
2, & q\equiv 3 \mod 4 \text{ and } k=9,11,17,\\
3, & q\equiv 3 \mod 4 \text{ and } k=7.
\end{cases}
\]
Define $n_+=22-\phi(k)-n_-$. Then
\[
R_a(T) = (1-q\,T)^{n_+}\,(1+q\,T)^{n_-}.
\]}
\end{Corollary}

\medskip

{\sl Proof:}
In the first cases, $\sqrt{-1}\in\F_q$ or all generators of $S_X$ can be defined over $\Q$. Hence Frob$_q^*$ acts as multiplication by $q$ on $S_X$. 

In all other cases, Gal$(\F_q(\sqrt{-1})/\F_q)$ acts nontrivially on the section, components of the fibres of type $IV$ or $IV^*$ resp.~on the lines $\ell_{\pm}$. The dimension $n_-$ of the $(-1)$-space of the Galois action is easily verified.
 \qed

\medskip

For $q\equiv 1 \mod m$, we combine Lemma \ref{Lem:R_t} (or~Remark \ref{Rem:k=3}) and Corollary \ref{Cor:NS} to obtain the $\zeta$-function of all K3 surfaces in Theorem 1 and 2. In particular, we deduce their modularity.

%This K3 surface $X$ is defined by
%$$x_0^6+x_0x_1^5+x_1x_2^5+x_3^2=0\in\P^3(1,1,1,3)$$
%of degree $6$ in the weighted projective $3$-space $\P^3(1,1,1,3)$.
%It is covered by the Delsarte surface
%$$V_6^2: x_0^6+x_0x_1^5+x_1x_2^5+x_3^6=0\subset\P^3$$
%with a covering map $(x_0,x_1,x_2,x_3)\mapsto (x_0,x_1,x_2,x_3^3)$.
%By Shioda \cite{Sh}, $V_6^2$ has the Picard rank $6$ with $NS(V_6^2)$
%spanned by $6$ lines:$L_{\zeta}: x_1=0,\, x_0-\zeta\,x_3=0$ where
%$\zeta^6=-1$. 
%Now we consider our K3 surface $X$.  The covering map sends these $6$ lines to
%$L_{\zeta}: x_1=0,\, x_0-\zeta\,x_3=$ with $(\zeta^2)^3=-1$. Therefore,
%$X$ has Picard number $2$.

\end{section}

\begin{section}{Galois representations}

In this section we will study the Galois representations of dimension $\phi(k)$ associated
to the transcendental parts of these K3 surfaces.
\medskip

First we recall the result of Nikulin \cite{Ni}. Most of the following statements hold true in full generality, but we only state the special case of minimal rank of the transcendental lattice which is relevant to our issues.

\begin{Theorem} [Nikulin] {\sl Let $X$ be an algebraic K3 surface. Then
$H_X$ is a finite cyclic group of order $k$. Assume
that $\phi(k)=\mbox{rank}(T_X)$. Thus $k$ takes the values listed in
Theorem 1 and Theorem 2.

(a) For those values of $k$, the ring $\Z[\zeta_k]$ is a PID, i.e. is of class
number one.  (See \cite{MM}, and also \cite{MO}.)

(b) The representation of $H_X=\langle g_k \rangle$ in $T_X\otimes\Q$ is isomorphic to
a direct sum of irreducible representations of $H_X$ of dimension one.
%(i.e., $\phi(k)=\mbox{rank}(T_X)$).

(c) Let $\Phi_k(x)$ denote the $k$-th cyclotomic polynomial. Regard $T_X\otimes \Q$ 
as a $\Z[\langle g_k\rangle]$-module via the natural action of $g_k$ on $T_X$.
Then $T_X\otimes \Q$ is a torsion free $\Z[\langle g_k\rangle ]/\langle\Phi_k(g_k)\rangle$-module. 
Identifying $\Z[\langle g_k\rangle]/\langle\Phi_k(g_k)\rangle$ with the ring of integers $\Z[\zeta_k]$, we derive an isomorphism 
$$T_X\otimes \Q\simeq \Q[\zeta_k]\quad\mbox{as $\Z[\zeta_k]$-modules}.$$
}
\end{Theorem}
\medskip

We also recall the result of Zarhin on the Hodge groups of complex K3 surfaces.
Recall that $T_X$  is the orthogonal
complement of $S_X$ with respect to the cup product which we denote by $<\cdot,\cdot>$.

\begin{Theorem} [Zarhin] 
{\sl Let $X$ be a complex K3 surface, and 
$\rm{Hdg}\subset \rm{Aut}(H^2(X,\mathbf{Q}))$ be the Hodge group. 
%Let $T_X$ be
%the group of transcendental cycles on $X$, that is, it is the orthogonal
%complement of $NS(X)$ with respect to the cup product $<\quad,\quad>$,
Then

(a) $T_X$ is an irreducible $\rm{Hdg}$-module,

(b) $E:=\rm{End}_{\rm{Hdg}}(T_X)$ is a commutative field. $E$ has an involution induced
by $<\cdot,\cdot>$ with totally real fixed field $E_0$, and either $E=E_0$ or 
$E$ is a totally complex quadratic extension of $E_0$, and $<\cdot,\cdot>$ induces
a symmetric, respectively, Hermitian form $\Phi: T_X\times T_X\to E$,

(c) $\rm{Hdg}= SO(T_X,\Phi)$, respectively, $U(T_X,\Phi)$,

(d) Let $U(X)$ be the image of $\rm{Aut}(X)$ in $\rm{Aut}(T_X)$. Then $U(X)$ is
contained in the roots of unity of $E$, and hence is cyclic of order $n$ for some $n$,
and $\phi(n)|[E:\Q]|\rm{dim}_{\Q}(T_X)$.}
\end{Theorem}

We will elaborate this theorem for our examples of K3 surfaces.

\begin{Corollary} {\sl Let $X$ be a K3 surface in Theorem 1 or Theorem 2.
Then $[E:\Q]=\phi(k)=\rm{dim}_{\Q}(T_X)$, so that $E$
is a cyclotomic field, and hence a CM field, over $\Q$
of degree $\phi(k)$.}
\end{Corollary}

This gives a Hodge theoretic proof that our K3 surfaces in Theorem 1 and Theorem 2
are all of CM type.  (Confer Corollary 1.) 
\medskip

Now we will study the Galois representations associated to our
K3 surfaces. The main point is the comparison of the Hodge structure
and some piece of the Galois representation.

\begin{Lemma} {\sl Let $X$ be one of the K3 surfaces in Theorem 1
and Theorem 2.  Then for each $k$, the Galois representation associated to $T_X$ has
dimension $\phi(k)$, and is irreducible over $\Q$. In fact this $\Gal(\Qbar/\Q)$
representation is induced
from the Jacobi sum Grossencharacters corresponding to the first character associated
to $X$ in the table of Section 5. }
\end{Lemma}

{\sl Proof}: 
Since $H^2(X)$ is a submotive of $H^2(\mathcal{F}_m)$, the results of Section 5,
in particular Remarks \ref{Rem7} and \ref{Rem8},
show that the $\Gal(\Qbar/\Q(\zeta_m))$ representation defined by $T_X$ is a
sum of one dimensional representations which are simply transitively permuted by 
$\Gal(\Q(\zeta_m)/\Q)$. The claim follows
\qed

\begin{Remark} 
{\rm
Let $K$ be a finite field extension of $\Q$. Since the $\Gal(\Qbar/\Q)$
representation given by $T_X$ is induced from a character (of
$\Gal(\Qbar/\Q(\zeta_m))$\/), it follows from Mackey's formula 
that the restriction of this representation to $\Gal(\Qbar/K)$ is a
direct sum of cyclic representations. In particular the $L$-function 
of the motive $T_X\otimes K$ and with it the $L$-function of $H^2(X\otimes K)$ 
is modular for any $K$.}
\end{Remark}

\begin{Remark}  
{\rm 
A related modularity result is
obtained by invoking the automorphic induction of Arthur and
Clozel \cite{AC90}.  
The first result we proved is the $GL(1)$-modularity for
a CM motive over $\Q(\zeta_k)$. 
The automorphic induction of Arthur and Clozel 
then takes this $GL(1)$ automorphic representation to
a $GL(\varphi(k))$ automorphic cuspidal representation over $\Q$ having the same
$L$-function (for the standard representation).
 
}
\end{Remark}

\end{section}

\begin{section}{Arithmetic mirror symmetry}

In this section, we will show that for all values except $k=3$, a mirror
K3 surfaces exist. For some of them, we interpret mirror symmetry arithmetically.

\medskip

In the literature, there are several variants of mirror symmetry for
K3 surfaces. Here we employ the notion of mirror symmetry for
lattice polarized K3 surfaces introduced by Dolgachev \cite{Dol}, based on
the Arnold strange duality.

\begin{Definition}
{\sl 
Let $X$ be an algebraic K3 surface. Then a K3 surface $\breve X$ is called a {\em mirror} of $X$ if
\begin{eqnarray}\label{eq:mirror}
T_X = U_2 \oplus S_{\breve X}.
\end{eqnarray}}
\end{Definition}

\medskip

Usually mirror symmetry is exhibited on the level of families of K3 surfaces. Here we are only interested in the existence of a mirror K3 surface. Then we will investigate arithmetic properties.

%

%For pseudo-ample lattice polarized K3 surfaces, the mirror symmetry
%is formulated as follows.

%Let $\mathcal{L}=U_2^3\oplus(-E_8)^2$
%be the K3 lattice.  Let $\mathcal{M}\subset\mathcal{L}$ be a primitive sublattice
%which has a unique embedding into $\mathcal{L}$ and let $\mathcal{M}^{\perp}$ be
%the orthogonal complement of $\mathcal{M}$ in $\mathcal{L}$. 
%Suppose that $\mathcal{N}\subset \mathcal{M}^{\perp}$ is a sublattice which
%admits an orthogonal splitting
%$$\mathcal{M}^{\perp}=U_2\oplus\mathcal{N}.$$
%Then there exists a pair of pseudo-ample lattice polarized K3 surfaces
%$X_{\mathcal{M}}$ and $X_{\mathcal{N}}$ such that
%$\mathcal{M}=NS(X_{\mathcal{M}})$ and $\mathcal{N}=NS(X_{\mathcal{N}})$ and that
%$$\rho(X_{\mathcal{M}})+ \rho(X_{\mathcal{N}})=20.$$

%Here is our results.

\begin{Lemma} 
{\sl 
For each of the K3 surfaces listed in Theorem 1 and 2 except for $k=3$, a
mirror K3 surface exists.}
\end{Lemma}

{\sl Proof}: We give an abstract proof, based on the following fact: If $T_X$ admits an orthogonal splitting
\begin{eqnarray}\label{eq:T}
T_X = U_2 + M
\end{eqnarray}
as in (\ref{eq:mirror}), then $M$ embeds primitively into $H^2(X,\Z)$. Since $M$ has signature $(1,19-\rho(X))$, there is a K3 surface $\breve X$ with $S_{\breve X}=M$ by the work of Nikulin \cite{Ni}. By definition, $\breve X$ is a mirror of $X$. Hence we have to show that for each K3 surface in Theorem 1 and 2 except for $k=3$, $T_X$ admits an orthogonal splitting (\ref{eq:T}).

\medskip

For each K3 surface in Theorem 1 and 2 except for $k= 3, 11, 19$, an orthogonal splitting (\ref{eq:T}) is given in the tables in section \ref{s:results}. For $k=3$, no such splitting exists since $T_X$ is positive definite. For $k=11, 19$, we use a primitive embedding 
\[
A_2 \hookrightarrow E_8.
\]
In the present cases, this induces a primitive embedding
\[
\langle 2 \rangle \oplus (-A_2) \hookrightarrow T_X.
\]
Hence the claim follows from the fact that 
\[
U_2 \hookrightarrow \langle 2 \rangle \oplus (-A_2).
\]
To see this, choose generators $b$ of the left summand and $a_1, a_2$ for the right summand with intersection form $\begin{pmatrix} -2 & 1\\1 & -2\end{pmatrix}$. Then $U_2$ can be identified with the sublattice $\langle b+a_1, b+a_1+a_2\rangle$. \qed

\begin{Remark}
{\rm
For each surface $X$ in Theorem 1 and 2 except for $k=3$, the mirror construction is in fact symmetric, i.e.~$X$ is a mirror of $\breve X$. This follows from lattice theory as used in section \ref{s:cycles}, since there is only one class per genus.}
\end{Remark}

For many surfaces, we can write down a mirror explicitly. For instance, for the K3 surfaces in Theorem 1, we have the following mirror pairs:
\begin{eqnarray*}
\{k=12\} & \leftrightarrow & \{k=44, 66\}\\
\{k=28, 36, 42\} & \leftrightarrow & \{k=28, 36, 42\}
\end{eqnarray*}
%K3 surfaces corresponding to $k=28, 36$ and $42$ are self-mirrors.
\medskip

The same applies to some of the K3 surface in Theorem 2:
\[
k=9 \leftrightarrow k=27,\;\;\; k=5 \leftrightarrow k=25,\;\;\; k=13 \leftrightarrow k=13.
\]
%Here K3 surfaces
These relations give rise to an arithmetical mirror symmetry:

\begin{Lemma}
For each K3 surface in Theorem 1 and 2 except for $k=3, 7, 11, 17, 19$, there is a mirror which is a Delsarte surface. In particular, this mirror K3 surface is again of CM-type and modular.
\end{Lemma}

\medskip

No arithmetical mirror symmetry seems to be known for the remaining K3 surfaces. Nonetheless we can look for mirrors explicitly (except for $k=3$). Here we briefly comment on this for two cases:

\medskip

For $k=7, 17$, mirror symmetry is particularly easy, since any mirror $\breve X$ admits an elliptic fibration with section due to the embedding
\[
U_2 \hookrightarrow S_{\breve X}.
\]
For $k=7$, e.g., we require exactly one reducible fibre with root lattice $A_1$ (i.e.~type~$I_2$ or $III$) and a section $P$ meeting this fibre in the non-zero component, but not intersecting the zero section. These conditions give rise to a 16-dimensional family of K3 surfaces as follows. 

The $A_1$ fibre (located at $\infty$) is encoded in the general Weierstrass equation with polynomial coefficients $b(t),c(t)\in K[t]$
\begin{eqnarray}\label{eq:fam}
y^2 = x^3 + a\,t^4\, x^2 + b(t)\, x + c(t),\;\;\; a\in K,\;\, \deg(b)\leq 7,\,\;\deg(c)\leq 10.
\end{eqnarray}
We require a  section  meeting the non-zero component of the $A_1$ fiber at $\infty$, but not intersecting the zero section.
Any such section $P$ can be given polynomially as
\begin{eqnarray}\label{eq:P}
P = (X(t), Y(t)),\;\;\; \deg(X(t))\leq 3, \;\deg(Y(t))\leq 5.
\end{eqnarray}
By inspection, the polynomials $X(t), Y(t), a, b(t)$ determine the remaining coefficient $c(t)$ of the Weierstrass form uniquely.
Taking the three normalisations due to M\"obius transformations and rescaling $x, y$ into account, we derive a 16-dimensional family of elliptic K3 surfaces as claimed.
Any general member of this family serves as a mirror of $X$. 

However, at this time we have no knowledge about the $\zeta$-functions of K3 surfaces of Picard number $\rho(X)=4$ unless they are Delsarte surfaces. In the above family, there is only one member with exactly one reducible fiber (of type $A_1$) which is a Delsarte surface:
\[
W:\;\;\; y^2 = x^3 + t^7\,x + 1.
\]
Here the section $P$ is given by $(0, 1)$. Nonetheless, $W$ is not a mirror of $X$ since it has Picard number $\rho(W)=16$. This can be proved along the lines of section \ref{s:motif} by interpreting $W$ as a quotient of the Fermat surface $\mathcal{F}_m$ of degree $m=42$. Hence no arithmetical mirror of $X$ is known.

\medskip

Similarly for $k=17$, we can write down a mirror family of elliptic K3 surfaces with trivial lattice
\[
U_2 \oplus (-A_1) \oplus (-A_2) \oplus (-E_8)
\]
and a section $P$ of height $h(P)=17/6$, but then there is no non-degenerate Delsarte surface in this family at all. 

\end{section}

\vskip 1cm

%\vfill
%\pagebreak

\centerline{\bf Acknowledgments}

During the preparation of this paper, 
M.~Sch\"utt held positions at Harvard University and University of Copenhagen.
N. Yui was a visiting researcher at IHES France (January--March 2007),
Hebrew University of Jerusalem (March 2007), DPMMS University of Cambridge 
(April--June 2007), and Tsuda College Japan (July--August 2008).  
We are grateful for the hospitality and support of each institution. 

\vskip 1cm

\bibliography{plain}

%\bibliographystyle{plain}
%\bibliography{CYMonodromy}

\end{document}